\documentclass[11pt]{article}
\usepackage{amsmath}
\usepackage{amssymb}
\usepackage{amsthm}
\usepackage{epsf}

\newcounter{AbcT}

\newtheorem {Theorem}    {Theorem}[section]
\newtheorem {Definition} {Definition} [section]
\newtheorem {Problem}    {Problem}
\newtheorem {Conjecture} {Conjecture}

\newtheorem {Lemma}      [Theorem]    {Lemma}

\newtheorem {Proposition}[Theorem]    {Proposition}

\newtheorem {Example}    [Theorem]    {Example}

\newcommand {\Heads}[1]   {\smallskip\pagebreak[1]\noindent{\bf #1{\hskip 0.2cm}
}}
\newcommand {\Head}[1]    {\Heads{#1:}}

\newcommand {\proofs}     {\proof}
\newcommand {\prooft}[1]  {\Head{Proof #1}\nopagebreak[2]}

\newcommand  {\QED}    {\def\qedsymbol{$\blacksquare$}\qed}

\newcommand{\Z}{{\bf{Z}}}
\newcommand{\ignore}[1]{}
\newcommand{\br}{{\mbox{{\rm br}}}}
\newcommand{\rank}{{\mbox{{\rm rank}}}}

\def\pho{\rho}
\def\eps{\epsilon}
\def\E{{\bf{E}}}
\def\P{{\bf{P}}}

\def\N{{\bf N}}
\def\M{{\bf{M}}}
\def\R{\hbox{I\kern-.2em\hbox{R}}}
\def\A{{\cal{A}}}

\def\F{{\cal{F}}}

\def\|{\, | \, }
\def\v0{{\bf 0}}
\def\one{{\bf 1}}
\def\0{\hat{0}}
\def\1{\hat{1}}

\def\pa{{\tt path}}
\def\path{{\tt path}}

\def\lam{\lambda}

\def\phi{\varphi}

\def\be{\begin{equation}}
\def\ee{\end{equation}}
\def\Reff{{\cal R}_{\rm eff}}
\def\lra{\leftrightarrow}
\def\one{{\bf 1}}
\def\NC{{\cal{N}}}

\begin{document}
\title{Information flow on trees}
\author{Elchanan Mossel \\ INRIA Rocquencourt and Microsoft Research \\ \\ Yuval Peres \\ Hebrew University and U.C. Berkeley}
\maketitle

\begin{abstract}

Consider a tree network $T$, where each edge acts as an independent copy of
a given channel $M$, and information is propagated from the root.
For which $T$ and $M$ does the configuration
 obtained at level $n$ of $T$
typically contain significant information on the root variable?
This problem arose independently in biology, information theory and
statistical physics.
\begin{itemize}
\item
For all $b$, we construct a channel
for which the variable at the root of the $b$-ary tree
is independent of the  configuration at level $2$ of that tree,
yet for sufficiently large $B>b$, the  mutual information between
the configuration at level $n$ of
the $B$-ary tree and the root variable is bounded away from zero.
This is related to certain secret-sharing protocols.

\item
We improve the upper bounds on information flow
for asymmetric binary channels (which correspond to the  Ising model
with an external field) and for symmetric $q$-ary channels
(which correspond to Potts models).

\item
Let $\lam_2(M)$ denote the second largest
eigenvalue of $M$, in absolute value. A CLT of Kesten and
Stigum~(1966) implies that if
$b |\lam_2(M)|^2 >1$, then the {\em census}
of the  variables at any level of
the $b$-ary tree, contains significant information on the root variable.
We establish a converse: if $b |\lam_2(M)|^2 < 1$, then the
census of the variables at level $n$ of the $b$-ary tree is
 asymptotically independent of the root variable.
This contrasts with examples where  $b |\lam_2(M)|^2 <1$,
yet the  {\em  configuration} at level $n$
is not asymptotically independent of the root variable.

\end{itemize}
\end{abstract}

\section{Introduction}

Consider a process in which information flows from the
root of a tree $T$ to other nodes of $T$.
Each edge of the tree acts as a channel
on a finite alphabet $\A=\{1,\ldots,k\}$.
Denote by $\M_{i,j}$ the transition probability from $i$ to $j$, and
by $M$ the random function (or channel) which satisfies for all $i$
and $j$ that $\P[M(i) = j] = \M_{i,j}$.
Let $\lam_2(M)$ denote the eigenvalue of $\M$ which has the
second largest absolute value ($\lam_2(M)$ may be negative or non-real).
At the root $\rho$ one of the symbols of $\A$ is chosen according to
some initial distribution. 
We denote this (random) symbol by $\sigma_{\rho}$.
This symbol is then propagated in the tree
as follows. For each vertex $v$ having as a parent $v'$, we let
$\sigma_v = M_{v',v}(\sigma_{v'})$, where the $\{M_{v',v}\}$
are independent copies of $M$.
Equivalently, for a vertex $v$, let $v'$ be the parent of $v$, and let
$\Gamma(v)$ be the set of all vertices which are connected to $\pho$
through paths which do not contain $v$. Then the process
satisfies:
\[
\P[\sigma_v = j | (\sigma_w)_{w \in \Gamma(v)}] = \P[\sigma_v = j | \sigma_{v'}]
 = \M_{\sigma_{v'},j}.
\]

It is very natural to study this process in the context of biology,
statistical physics and communication theory. See
\cite{EKPS},\cite{M:lam2} and the references there for more background.

Let $d(,)$ denote the graph-metric distance on $T$, and
$L_n = \{v \in V: d(\rho,v) = n\}$ be the $n$'th level of the tree.
For $v \in V$ and $e=(v,w) \in E$ we denote $|v| = d(\rho,v)$ and
$|e| = \max\{|v|,|w|\}$.
We denote by $\sigma_n = (\sigma(v))_{v \in L_n}$ the symbols at the
$n$'th level of the tree.
We let $c_n = (c_n(1),\ldots,c_n(k))$ where
\[
c_n(i) = \#\{v \in L_n : \sigma(v) = i\}.
\]
In other words, $c_n$ is the {\sl census} of the $n$'th level.
Note that both $(\sigma_n)_{n=1}^{\infty}$ and
$(c_n)_{n=1}^{\infty}$ are markov chains.

For distributions $P$ and $Q$ on the same space the total
variation distance between $P$ and $Q$ is
\begin{equation} \label{eq:totalvar}
D_V(P,Q) = \frac{1}{2} \sum_{\sigma} |P(\sigma) - Q(\sigma)|.
\end{equation}

\begin{Definition} \label{def:reconstruction}
The reconstruction problem for $T$ and $M$ is {\bf solvable} if
there exist $i,j \in \A$ for which
\begin{equation} \label{cond:arb_l1}
\lim_{n \to \infty} D_V(\P_n^i,\P_n^j) > 0,
\end{equation}
where
$\P_n^\ell$ denotes the conditional distribution of
$\sigma_n$ given that $\sigma_{\rho} = \ell$.
\end{Definition}

\begin{Definition} \label{def:count_reconstruction}
The reconstruction problem for $T$ and $M$ is {\bf census-solvable} if
there exist $i,j \in \A$ for which
\begin{equation} \label{cond:arb_count_l1}
\lim_{n \to \infty} D_V(\P^{(c),i}_n, \P^{(c),j}_n) > 0,
\end{equation}
where
$\P^{(c),\ell}_n$ denotes the conditional distribution of
$c_n$ given that $\sigma_{\rho} = \ell$.
\end{Definition}

Assume that $\P[\sigma_\rho=i]>0$ for every $i \in \A$.
Then the following conditions are equivalent, see e.g. \cite{M:lam2}:
\begin{itemize}
\item
The reconstruction problem for $T$ and $M$ is not solvable.
\item
The configurations
$\sigma_n$  are asymptotically independent of $\sigma_\rho$, i.e.,
$\lim_{n \to \infty} I(\sigma_\rho,\sigma_n)=0$, where $I$ is the mutual
information operator.
\item
The sequence $(\sigma_n)_{n=1}^{\infty}$ has a trivial tail $\sigma$-field.
\item
The distribution of $\{\sigma_v\}_{v \in T}$
is an extremal Gibbs measure.
\end{itemize}

Similar conditions apply to census-solvability, where $\{c_n\}$
replace $\{\sigma_n\}$.
Note that if the reconstruction problem is census-solvable
it is also solvable.

The reconstruction problem was first studied in statistical physics
for the Ising model on the tree, in the equivalent form as a question
regarding the extremality of the free Gibbs measure for this model.
See \cite{Sp},\cite{Hi}, \cite{BRZ},\cite{I1}, \cite{I2}
and \cite{KMP}. Reference \cite{EKPS} contains extensive background on
this problem for the Ising model.
In these papers the reconstruction problem for the Ising model on trees
was solved. Writing $\lam_2(M)$ for the second eigenvalue
of $\left( \begin{array}{ll} 1 - \eps & \eps \\
\eps & 1 - \eps \end{array} \right)$ (that is,
$\lam_2(M) = 1 - 2 \eps$), the reconstruction problem
is solvable for the $b$-ary tree $T_b=(V_b,E_b)$ if and only if
$b \lam_2^2(M) > 1$. An analogous threshold
for general trees was established
in \cite{EKPS}.

A few combinatorial questions arise naturally
in the context of these tree processes
(see also \cite{BW1}, \cite{BW2} and \cite{BW3}).
In Section \ref{sec:states_and_trees} we discuss the following question.
For a channel $M$ we ask does there exist any $b$
such that the reconstruction problem is solvable for the infinite
$b$-ary tree $T_b$ and the channel $M$?
In Theorem \ref{thm:dist_comb} we present a criteria for deciding
this problem for a channel $M$. Using this criteria we construct a channel
with the following properties:
\begin{Example} \label{ex:strange_shift}
Write $M^h$ for the $h$'th iteration of the channel $M$.
Then for all $h > 1$ there exists a channel $M$ such that
\begin{itemize}
\item $M^h(j)$ has the same distribution for all $j \in A$.
\item For all $h' < h$ there exist $i$ and $j$ such that $M^{h'}(i)$
  and $M^{h'}(j)$ have different distributions.
\item
When $b$ is large the reconstruction problem is solvable for the tree
$T_b$
and the channel $M$.
\end{itemize}
\end{Example}
This is a generalization of an example which appeared
in \cite{M:recur} and \cite{M:lam2}.
In Section \ref{sec:states_and_trees2} we give more delicate
constructions which are related to secret-sharing protocols (see \cite{Sh}).
\begin{Theorem} \label{thm:crypt}
Let $b > 1$ be an integer and $T$ be the $2$-level $b$-ary tree.
There exist a channel $M$ such that for any initial distribution,
$\sigma_r$ and $\sigma_{\partial}$ are independent
(where $\sigma_{\partial}$ is the
configuration at the leaves of the $2$-level $b$-ary tree), yet
when $B$ is sufficiently large, the reconstruction problem for the channel $M$
and the $B$-ary tree $T_B$ is solvable.
\end{Theorem}

It is tempting to try to find thresholds for the reconstruction problem which
depend only on $b$ and $\lam_2(M)$.
For binary symmetric channels the threshold for reconstruction is
$b \lam_2^2(M) = 1$
(this is also the threshold for census-reconstruction for general channels:
see theorems \ref{thm:blam2} and \ref{thm:non_count_rec} below).

In \cite{M:lam2} it is shown that for some natural generalizations of the binary symmetric
channel, the threshold $b \lam_2^2(M) = 1$ is not the threshold for reconstruction.
In particular, it is shown that for the Potts model (the $q$-ary symmetric channel),
where the transition matrix is:
\begin{equation} \label{eq:potts1}
\M =
 \left(
\begin{array}{lllll} 1 - (q-1) \delta & \delta & \ldots     & \delta   \\
                    \delta & 1 - (q-1) \delta &  \delta       & \ldots \\
                    \vdots & \ldots & \ddots & \vdots            \\
                    \delta & \ldots & \delta     & 1 - (q-1) \delta \end{array}
        \right),
\end{equation}
if $b \lam_2(M) = b(1 - q \delta) > 1$ and $q$ is sufficiently large, then the reconstruction
problem is solvable. Similar results hold for asymmetric binary channels.
On the other hand, it is well known, see e.g. \cite{M:lam2},
that for both families
the reconstruction problem is unsolvable if $0 \leq b \lam_2(M) \leq 1$.
In Section
\ref{sec:improved} we improve the bound $b \lam_2(M) \leq 1$ by proving:
\begin{Proposition} \label{prop:potts_no_reconstruct}
Let $\M$ be the matrix (\ref{eq:potts1}).
Then the reconstruction problem for $T_b$ and $M$ is unsolvable when
\begin{equation}
b \frac{(1 - q \delta)^2}{1 - (q-2) \delta} \leq 1.
\end{equation}
\end{Proposition}
Similar results hold for asymmetric binary channels.
The proofs are based on a reduction to symmetric binary channels on general
trees which were analyzed in \cite{EKPS}.

In sections \ref{sec:general_count}, \ref{sec:easy_count} and
\ref{sec:potts_count} we consider census reconstruction and show that the threshold
for census reconstruction is $b |\lam_2(M)|^2 = 1$ .
From \cite{KS} it follows that if $b |\lam_2(M)|^2 > 1$,
then the reconstruction problem is census solvable.
We extend this result to general trees by proving:
\begin{Theorem} \label{thm:blam2}
Let $T$ be an infinite tree and write $\br(T)$ for the branching number of the tree.
Let $M$ be a channel such that
$\br(T) |\lam_2(M)|^2 > 1$, then the reconstruction problem is solvable
for $T$ and $M$.
\end{Theorem}
The results of \cite{KS} also plays a crucial role in proving:
\begin{Theorem} \label{thm:non_count_rec}
The reconstruction problem for the $b$-ary tree and the channel $M$ is not census-solvable
if $b |\lam_2(M)|^2 < 1$.
\end{Theorem}
We conjecture that the result of Theorem \ref{thm:non_count_rec}
should hold also when $b |\lam_2(M)|^2 = 1$, but we have verified this only for
the following channels:
\begin{Theorem} \label{thm:no_count_reconstruct}
Let $M$ be the $q$ state Potts model, or the asymmetric Ising model,
and suppose that $b \lam_2^2(M) \leq 1$, then the reconstruction problem
is unsolvable for the $b$-ary tree and the channel $M$.
\end{Theorem}

In Section \ref{sec:open} we discuss some open problems.

Some of the results in this paper concern general trees.
For an infinite tree $T$
many of its probabilistic properties are determined by the branching
number $\br(T)$.
This is the supremum of the real numbers  $\lam \ge 1$,
such that $T$ admits a positive flow from the root to infinity, where
on every edge $e$ of $T$, the flow is bounded by $\lam^{-|e|}$.
Here $|e|$ denotes  the number of edges, (including $e$) on the path
from $e$ to the root; $\br(T)^{-1}$ is the critical probability
for Bernoulli percolation on $T$.
See \cite{Lyons} and \cite{EKPS} for equivalent definitions
of $\br(T)$  in terms of percolation, cutset sums and electrical conductance.
We note that for the regular tree $T_b$ we have $\br(T_b) = b$.

\section{Distinguishing states by tree networks} \label{sec:states_and_trees}

Let $M$ be a channel on an alphabet $\A$ of size $k$ with $\lam_2(M) = 0$.
Looking at the Jordan form of $\M$ we see that $\rank(\M^k) = 1$, and therefore
$M^k(i)$ has the same distribution for all $i \in \A$.
Moreover, by Theorem \ref{thm:non_count_rec} it follows that for all $b$, the reconstruction
problem is not {\em census} solvable for the channel $M$ and the tree $T_b$.
Does this mean that the reconstruction problem is unsolvable for $T_b$ and
the channel $M$? In this section we answer this question.

\begin{Definition} \label{def:dist_comb}
Let $M$ be a channel on $\A$ and consider the minimal
equivalence relation $\sim$ on $\A$ which satisfies:
If $i$ and $j$ satisfy for all $\ell$,
\begin{equation}\label{eq:comb_cond}
\sum_{\ell' \sim \ell} \M_{i,\ell'} = \sum_{\ell' \sim \ell} \M_{j,\ell'},
\end{equation}
then $i \sim j$.
If $i \sim j$, we say that $i$ and $j$ are {\em indistinguishable}.
Otherwise, we say that $i$ and $j$ are {\em distinguishable}.
\end{Definition}

The motivation for this definition is that if
$\M_{i,\ell} = \M_{j,\ell}$ for all $\ell$ , then after one application
of the random function $M$,
one cannot distinguish
between the hypotheses that $i$ was the input and
the hypotheses that $j$ was the input.

\begin{Theorem} \label{thm:dist_comb}
If the states $i$ and $j$ are indistinguishable,
then there exists $N$ such that
for all $n \geq N$ and for all infinite rooted trees $T$:
\begin{equation} \label{eq:comb_indistinguishable}
D_V(\P_n^i,\P_n^j) = 0.
\end{equation}
On the other hand, for every channel $M$,
there exits $b$, such that for the tree $T_b$:
\begin{equation} \label{eq:comb_distinguishable}
\inf_{n \geq 1, i \nsim j} D_V(\P_n^i,\P_n^j) > 0.
\end{equation}
\end{Theorem}

\begin{Example} \label{ex:shift}
Let $\{Z_i\}_{i=1}^{\infty}$ be an i.i.d sequence of variables such that
$\P[Z_i = 0] = 1 - P[Z_i = 1] = p$ where $0 < p < 1$.
Let $h \geq 1$ and consider the channel $M$ defined by the markov chain
$Y_i = (Z_i,\ldots,Z_{i+h})$.
Thus $M$ has state space $\{0,1\}^{h+1}$ with the product $(p,1-p)$
probability measure.
It is easily seen that all the states of $M$ are
indistinguishable. Therefore by Theorem \ref{thm:dist_comb} there exists $N$ such
that for all trees and all $n \geq N$, it holds that $D_V(\P_n^i,\P_n^j) = 0$
(one can take $N=h$).
\end{Example}

\medskip\noindent
{\bf Example \ref{ex:strange_shift}.}
{\it
Let $\{Z_i\}_{i=1}^{\infty}$ be as in Example \ref{ex:shift}.
Set
\[
Y_i = \max_{0 \leq j \leq h} \{Z_i = \cdots = Z_{i+j} = 1\}.
\]
It is easily seen that $Y_i$ defines a channel $M$ on the space
$\{0,\ldots,h\}$.
Moreover, it is clear that for all $\ell$ the variables $\{Z_i\}_{i \geq \ell+h+1}$ and $\{Y_i\}_{i \leq \ell}$
are independent.
Therefore $\{Y_i\}_{i \geq \ell+h+1}$ and $\{Y_i\}_{i \leq \ell}$ are independent.
It follows that the variables $M^{h+1}(j)$ have the same distribution for all $j$.
Thus $\rank(\M^{h+1}) = 1$, and $\lam_2(M) = 0$.
Writing $\M$:
\[
\M = \left( \begin{array}{lllll}    p &   p(1-p) &    p (1-p)^2 &     \ldots &    (1-p)^h \\
                                    1   &   0 &         \ldots &             \ldots &    0   \\
                                    0   &   1   &       0 &             \ldots &         \vdots \\
                                    \ldots & \ddots & \ddots &         \ddots &    \vdots \\
                                    0 &     \ldots  &  0     &          p &       1-p \\
            \end{array} \right),
\]
we see that all the states of $M$ are distinguishable.
Therefore by Theorem \ref{thm:dist_comb}, when $b$ is sufficiently large we have
for the tree $T_b$ that $\inf_{n \geq 1,i,j} D_V(\P_n^i,\P_n^j) > 0$.
Note that the channel above is obtained by lumping the channel of example
\ref{ex:shift}.
This is a generalization of a channel appearing in \cite{M:recur};
see also \cite{LP2}.
}

\begin{Remarks}
\begin{itemize}
\item
Theorem \ref{thm:dist_comb} implies that if
for a channel $M$ all
the states are indistinguishable, then $\lam_2(M) = 0$.
\item
Suppose that the channel $M$ is reversible and has $\lam_2(M) = 0$.
Then looking at the diagonal form of $\M$, it follows that $\rank(\M) = 1$.
In particular, all the states of $M$ are indistinguishable.
\end{itemize}
\end{Remarks}

\prooft{of Theorem \ref{thm:dist_comb}}
The first claim is easy. We define a new equivalence relation between
states in $\A$. We let $i$ and $j$ be equivalent if there exists an $N$
such that for all trees and all $n \geq N$ we have
$D_V(\P_n^i,\P_n^j) = 0$. It is clear that this is indeed an equivalence
relation and that it satisfies (\ref{eq:comb_cond}). The proof of the
first claim follows.

For the proof of the second claim
let $\A'$ be the set of equivalence classes for $\sim$.
Let $M'$ be the channel on $\A'$ defined as follows.
For $i',j' \in \A'$ we choose $i \in i'$
and define
\begin{equation} \label{eq:defM'}
\P[M'(i') = j'] = \sum_{j \in j'} \P[M(i) = j].
\end{equation}
It is clear that (\ref{eq:defM'}) does not depend on the choice of $i$.
Moreover, it is clear that for any tree $T$, we may couple
the $M$ tree process and the $M'$ tree process in such a way that
for all $v$ in $T$ we have that $\sigma_v \in \sigma_v'$.
In particular, it follows that for all trees and all
$i \in i'$ and $j \in j'$:
\[
D_V(\P^{i'}_n,\P^{j'}_n) \leq D_V(\P^i_n,\P^j_n).
\]
Therefore, in order to prove the second claim it suffices to show that
there exists a $b$ such that for the tree $T_b$ and the channel $M'$ we have,
\begin{equation} \label{eq:glue_dist}
\inf_{i,j \in \A', n} D_V(\P^i_n,\P^j_n) > 0.
\end{equation}
The crucial property of $\M'$ we will exploit is that for
all $i \neq j \in \A'$ we have $(\M_{i,\ell}')\neq (\M_{j,\ell}')$.
Thus in order to simplify the presentation of the proof, we will assume
that the channel $\M$ on the alphabet $\A$ satisfies
$(\M_{i,\ell})_{\ell = 1}^{|\A|} \neq (\M_{j,\ell})_{\ell = 1}^{|\A|}$
for all $i \neq j \in \A$

We are going to show that for every $\eps > 0$, there exits $b=b(\eps)$,
and a recursive algorithm, which given the symbols at the $n$'th
level of $T_b$, reconstruct the symbol at the root with probability at least
$1 - \eps$ (uniformly for all initial distributions of the root).
This implies (see \cite{M:lam2}) that (\ref{eq:glue_dist}) holds.

We let $\delta = \min_{i,j} \|(\M_{i,\ell}) - (\M_{j,\ell})\|_{\infty}$,
and assume that $\eps < \delta/8$.
By standard results in the theory of large deviations (see e.g. \cite{DZ}),
there exits a positive
constant $C$ such that for all $i$,
if $X_i$ is the empirical distribution vector
of $b$ independent trials with distribution $\M_{i,\ell}$, then,
\begin{equation} \label{eq:LLD1}
\P[\|X_i - b \M_{i,\ell}\|_{\infty}\| > b \delta/8] \leq
\exp(-C b \delta).
\end{equation}
Similarly, if $X$ is the sum of $b$ i.i.d. $\{0,1\}$
variables, each of which has the value $1$ with probability
$\eps$ and the value $0$ with probability $1-\eps$, then,
\begin{equation} \label{eq:LLD2}
\P[|X - \eps b| > b \delta/8] \leq \exp(-C b \delta).
\end{equation}
We now choose $b$ such that $\exp(-C b \delta) \leq \eps/2$,
and apply the following recursive algorithm in order to reconstruct the
symbol at the root of the tree:
\begin{itemize}
\item
For each vertex $v$, construct the empirical distribution vector
$X_v$ of the reconstructed values for the children of that vertex.
\item
Reconstruct at $v$ the value $i$, where $i$ minimizes the
distance $\|X_v - b \M_{i,\ell}\|_{\infty}$.
\end{itemize}
It now follows by induction, that the probability of correct reconstruction
is at least $1-\eps$. Indeed if this is true for all the children of $v$,
then by (\ref{eq:LLD2}), with probability at least $1-\eps/2$, for at least
$(1 - \delta/4)b$ of the children we reconstructed the correct value.
Thus by (\ref{eq:LLD1}), with probability at least $1-\eps$, the vector
$X_v$ satisfies $\|X_v - b \M_{i,\ell}\|_{\infty} < b \delta/2$ where $i$ is
the symbol at $v$. It now follows that we reconstruct the correct value with
probability at least $1 - \eps$.
\QED

\section{The distinguishing power of tree networks}
\label{sec:states_and_trees2}

Let $M$ be a channel on $\A$ and
let $T$ be a tree which consists of $2$ nodes, the root $r$
and an additional vertex $v$ (the $1$-ary $1$-level tree).
Suppose that $I(\sigma_r,\sigma_{\partial T}) = 0$, where $\partial T$ are the
vertices in the boundary of the tree $T$ (in this case $v$), and $I$
is the mutual information operator.
The assumption $I(\sigma_r,\sigma_{\partial T}) = 0$ is
equivalent to $\rank(\M) = 1$.
Therefore, all the states of $M$ are indistinguishable and
the reconstruction problem for $T_b$ and the channel $M$ is unsolvable for all $b$.

On the other hand, let $T$ be the $2$-level $1$-ary tree.
Example \ref{ex:strange_shift} is
a channel $M$ such that for the tree $T$ we have
$I(\sigma_r,\sigma_{\partial T}) = 0$, yet
the reconstruction problem for $T_b$
is solvable when $b$ is sufficiently large.

It is natural to ask similar question for other finite trees $T$:
Suppose that the channel $M$
satisfies for any initial distribution $I(\sigma_r,\sigma_{\partial T}) = 0$,
does this imply that the reconstruction problem for $T_b$ is unsolvable for all $b$?

Arguing as before, it is clear that this is the case if there exists
a vertex in $\partial T$ which is at distance $1$ from the root.
We show that these are the only trees for which the $0$-information condition
implies non-reconstruction for all $b$.

In the following theorem we give a construction
related to secret-sharing protocols \cite{Sh}.

\medskip\noindent
{\bf Theorem \ref{thm:crypt}}
{\it
Let $b > 1$ be an integer and $T$ be the $2$-level $b$-ary tree.
There exists a channel $M_b$ such that for any initial distribution
$I(\sigma_r,\sigma_{\partial T}) = 0$, yet
when $B$ is sufficiently large, the reconstruction problem for the channel $M$
and $T_B$ is solvable.
}

\begin{Definition} \label{def:crypt}
Let $\F$ be a finite field with $q > b+2$ elements.
Let $x_1,\ldots,x_{b+1}$ be a fixed set of non-zero element of $\F$.
We define a channel on the state space
\[
\F^b[x] = \{ f(x) : f(x) \in \F[x], \deg f \leq b \}.
\]
Given $f$, take $I$ to be a uniform variable in the set $\{1,\ldots,b+1\}$,
then take $M(f)$ to be $g \in \F^b[x]$ which satisfy $g(0) = f(x_I)$
with probability
\[
\left|\{g \in \F^b[x] : g(0) = f(x_I) \} \right|^{-1} = q^{-b}.
\]
\end{Definition}

\begin{Proposition} \label{lem:crypt_smoothness}
Let $\{x_1,\ldots,x_{b+1}\}$ be a set of elements of $\F$.
There is a one to one mapping from $\F^b[x]$ to $\F^{b+1}$
defined by
\begin{equation} \label{eq:map_poly_roots}
f \to \left( f(x_1),\ldots,f(x_{b+1}) \right).
\end{equation}
The inverse map is defined by the interpolation polynomial:
\begin{equation} \label{eq:map_poly_root_inv}
(y_1,\ldots,y_{b+1}) \to f(x) = \sum_{i=1}^{b+1}
\frac{\prod_{j \neq i} (x - x_j)}{\prod_{j \neq i} (x_i - x_j)} y_i.
\end{equation}
\end{Proposition}

\begin{Lemma} \label{lem:crypt_no_reconsturct}
Let $T$ be the $2$-level $b$-ary tree. Then for any initial distribution,
$I(\sigma_r,\sigma_{\partial}) = 0$.
\end{Lemma}

\begin{proof}
We show that for all $f \in \F^b[x]$ and
$h = (h_{i,j})_{i,j=1}^b \in (\F^b[x])^{b^2}$ it holds that
\[
\P[\sigma_{\partial T} = h | \sigma_r = f]
\]
is independent of $f$.
This would imply the claim of the Lemma.
We denote by $(v_i)_{i=1}^b$ the children of $r$,
and by $(w_{i,j})_{j=1}^b$ the children of $v_i$.
We then have
\begin{eqnarray} \label{eq:cond_rule_crypt}
\P[\sigma_{\partial T} = h | \sigma_r = f] &=&
\sum_{a_1,\ldots,a_b \in F} \P[\forall i, \sigma_{v_i}(0) = a_i | \sigma_r = f]
\P[\sigma_{\partial T} = h | \forall i, \sigma_{v_i}(0) = a_i]
\\ \nonumber &=&
\sum_{a_1,\ldots,a_b \in F} \P[\forall i, \sigma_{v_i}(0) = a_i | \sigma_r = f]
\prod_{j=1}^b \P[\forall j, \sigma_{w_{i,j}} = h_{i,j} | \sigma_{v_i}(0) = a_i].
\\ \nonumber &=&
q^{-b^2} \sum_{a_1,\ldots,a_b \in F} \P[\forall i, \sigma_{v_i}(0) = a_i | \sigma_r = f]
\prod_{j=1}^b \P[\forall j, \sigma_{w_{i,j}}(0) = h_{i,j}(0) | \sigma_{v_i}(0) = a_i].
\end{eqnarray}
Note that Lemma \ref{lem:crypt_smoothness} implies that
if $x_1,\ldots,x_j$ is
set of nonzero elements and $1 \leq j \leq b$, then the uniform distribution measure on $\F^b[x]$
satisfies for all $y_1,\ldots,y_j,a \in \F$ that
\[
\P[\forall \, 1 \leq i \leq j, f(x_i) = y_i | f(0) = a] = q^{-j}.
\]
Thus the probability at the right hand side of (\ref{eq:cond_rule_crypt}):
\[
\P[\forall j, \sigma_{w_{i,j}}(0) = h_{i,j}(0) | \sigma_{v_i}(0) = a_i].
\]
does not depend on $a_i$. Now the claim follows.
\end{proof}

\begin{Lemma} \label{lem:crypt_perm}
Let $\{x_1,\ldots,x_{b+1}\}$ be a fixed set of non-zero element in $\F$.
For an element $f \in \F^b[x]$ satisfying $f(x_k) = y_k$ for $1 \leq k \leq b+1$
and a permutation $\pi \in S_{b+1}$ we define $f_{\pi}$ to be the element of $F^b[x]$
which satisfy
\begin{equation} \label{eq:def_f_pi}
f_{\pi}(x_{\pi(k)}) = y_k.
\end{equation}
Moreover, for $a,b \in \F$ define:
\begin{equation} \label{eq:def_nab}
n_{a,b} = \left| \left\{(f,\pi) \in F^b[x] \times S_{b+1} : f(0) = a, f_{\pi}(0) = b \right\} \right|.
\end{equation}
We then have for all $a \neq b, c \neq d$ that
\begin{equation} \label{eq:trans_inv_nab}
\begin{array}{ll} n_{a,b} = n_{c,d}, & n_{a,a} = n_{c,c} \end{array}.
\end{equation}
Moreover, for all $a \neq b$,
\begin{equation} \label{eq:crypt_unsmoothness}
n_{a,a} > n_{a,b}.
\end{equation}
\end{Lemma}

\begin{proof}
We note that if $c \in F$, and $f \in F^b[x]$, then $(f+c)_{\pi} = f_{\pi} + c$.
Similarly, if $0 \neq d \in F$, then $(df)_{\pi} = d f_{\pi}$.
Thus, for all $a,b$ we have $n_{a+c,b+c} = n_{a,b}$ and $n_{da,db} = n_{a,b}$.
Relations (\ref{eq:trans_inv_nab}) now follows.

Looking at (\ref{eq:map_poly_root_inv}) one sees that the elements of $\F^b[x]$
which satisfy $g(0) = a$ are exactly those elements which satisfy the equation
\[
\sum_{i = 1}^{b+1} c_i g(x_i) = a,
\]
where $c_i$ are some none zero constants.

Thus the elements $g$ of $\F^b[x]$ which satisfy
$g(0) = a$ and $g_{\pi}(0) = b$
are exactly those elements of $\F^b[x]$ which
satisfy the equations
\begin{equation} \label{eq:lin_crypt2}
\sum_{i=1}^{b+1} c_i g(x_i) = a,
\end{equation}
and
\begin{equation} \label{eq:lin_crypt1}
\sum_{i=1}^{b+1} c_{\pi(i)} g(x_i) = b.
\end{equation}

Note that by Proposition \ref{lem:crypt_smoothness} the number of solutions of these
equations, $s_{\pi}(a,b)$, in $\F^b[x]$ is the same as the number of solutions of these equation over $\F$ in
the variables $g(x_1),\ldots,g(x_{b+1})$.
It now follows that for every permutation $\pi \neq 1$ and $a \neq b$:
\begin{equation} \label{eq:all_pi}
s_{\pi}(a,b) \leq s_{\pi}(a,a).
\end{equation}
Moreover, when $\pi = 1$,
\begin{equation} \label{eq:id_pi}
0 = s_{\pi}(a,b)  < s_{\pi}(a,a) = q^{b}.
\end{equation}
Now,
\[
n_{a,b} = \sum_{\pi \in S_{b+1}} s_{\pi}(a,b) < \sum_{\pi \in S_{b+1}} s_{\pi}(a,a) = n_{a,a},
\]
as needed.
\end{proof}.

\begin{Lemma} \label{lem:crypt_reconstruct}
Let $f,g \in \F^b[x]$. Then $f$ and $g$ are indistinguishable in the sense of
definition \ref{def:dist_comb} if and only if $f = g_{\pi}$ for some $\pi \in S_{b+1}$.
\end{Lemma}

\begin{proof}
We will write $f \sim_1 g$ to denote that $f$ and $g$ are indistinguishable,
and $f \sim_2 g$ when there exist $\pi \in S_{b+1}$ such that $f = g_{\pi}$.
It is clear that if $f \sim_2 g$, then $f \sim_1 g$.

On the other hand, suppose that $f \nsim_2 g$.
We will show that there exist an $h \in \F^b[x]$ such that
\begin{equation} \label{eq:c_pi_sum}
\sum_{h' \sim_2 h} \P[M(f) = h'] \neq \sum_{h' \sim_2 h} \P[M(g) = h].
\end{equation}
This would imply that indeed $\sim_1 = \sim_2$ (see definition \ref{def:dist_comb}).
We may write (\ref{eq:c_pi_sum}) equivalently,
\begin{equation} \label{eq:crypt_pi_sum}
\sum_{\pi \in S_{b+1}} \P[M(f) = h_{\pi}] \neq \sum_{\pi \in S_{b+1}} \P[M(g) = h_{\pi}].
\end{equation}
Writing $n_f(a) = |\{1 \leq j \leq b+1 : f(x_j) = a\}|$,
we may write the left hand side of (\ref{eq:crypt_pi_sum}) as:
\begin{equation} \label{eq:crypt_pi_sum2}
\frac{1}{(b+1) q^b} \sum_{\pi \in S_{b+1}} n_f(h_{\pi}(0)).
\end{equation}
Summing over all $h \in \F^b[x]$ with $h(0) = a$ in (\ref{eq:crypt_pi_sum2}) we obtain:
\begin{equation} \label{eq:crypt_pi_sum3}
\frac{n_{a,a} n_f(a) + \sum_{b \neq a} n_{a,b} n_f(b)}{(b+1) q^b} =
\frac{(n_{a,a} - n_{a,a+1}) n_f(a) + (b+1) n_{a,a+1}}{(b+1) q^b}.
\end{equation}
(by (\ref{eq:trans_inv_nab}) we may write $n_{a,a+1}$ for $n_{a,b}$ when $a \neq b$).
Since $f \nsim_2 g$ there exists $a \in \F$ such that $n_f(a) > n_g(a)$.
Thus:
\begin{eqnarray} \label{eq:crypt_pi_sum4}
\sum_{h : h(0) = a, \pi \in S_{b+1}} \P[M(f) = h_{\pi}] &=&
\frac{(n_{a,a} - n_{a,a+1}) n_f(a) + (b+1) n_{a,a+1}}{(b+1) q^b} \\ \nonumber &>&
\frac{(n_{a,a} - n_{a,a+1}) n_g(a) + (b+1) n_{a,a+1}}{(b+1) q^b} =
\sum_{h : h(0) = a, \pi \in S_{b+1}} \P[M(g) = h_{\pi}].
\end{eqnarray}
where we have used the fact that by (\ref{eq:crypt_unsmoothness}) $n_{a,a} > n_{a,a+1}$.
It now follows that there exists an $h$ for which (\ref{eq:crypt_pi_sum}) holds.
\end{proof}

\prooft{of Theorem \ref{thm:crypt}}
Fix $b$ and let $T$ be the $2$-level $b$-ary tree.
Let $M$ be the channel defined at definition \ref{def:crypt}.
Lemma \ref{lem:crypt_smoothness} implies that for any initial distribution $I(\sigma_r,\sigma_{\partial}) = 0$.

On the other hand by Lemma \ref{lem:crypt_reconstruct},
the channel $M$ satisfies that $f$ and $g$ are indistinguishable for $M$
if and only if $g = f_{\pi}$ for some $\pi$. It follows by Theorem \ref{thm:dist_comb} that when $B$ is sufficiently
large the reconstruction problem for $T_B$ is solvable (moreover, when $B$ is sufficiently large we may reconstruct
$\{f_{\pi}\}_{\pi \in S_{b+1}}$ the equivalence class of $f$, with probability close to $1$).
\QED
\section{Improved bounds for Potts models} \label{sec:improved}

In this section we will focus on the following two families of channels:
\begin{itemize}
\item
The asymmetric binary channels. These channels have the states space
$\{0,1\}$ and the matrices:
\begin{equation} \label{eq:assymetric}
\M = \left( \begin{array}{ll} 1 - \delta_1 & \delta_1
                           \\ 1 - \delta_2 & \delta_2
           \end{array} \right),
\end{equation}
with $\lam_2(M) = \delta_2 - \delta_1$.
\item
The symmetric channels on $q$ symbols. These have the state space
$\{1,\ldots,q\}$ and the matrices:
\begin{equation} \label{eq:potts}
\M =
 \left(
\begin{array}{lllll} 1 - (q-1) \delta & \delta & \ldots     & \delta   \\
                    \delta & 1 - (q-1) \delta &  \delta       & \ldots \\
                    \vdots & \ldots & \ddots & \vdots            \\
                    \delta & \ldots & \delta     & 1 - (q-1) \delta \end{array}
        \right),
\end{equation}
with $\lam_2(M) = 1 - q \delta$.
\end{itemize}
The reconstruction problem is unsolvable for (\ref{eq:assymetric}) when
$|b \lam_2(M)| \leq 1$ and unsolvable for (\ref{eq:potts}) when
$0 \leq b \lam_2(M) \leq 1$ (see Propositions 3 and 4 in \cite{M:lam2})
On the other hand,
Proposition \ref{prop:dlam2} implies that when $b \lam_2^2(M) > 1$ the
reconstruction problem is solvable for these channels.
Moreover, the main results of \cite{M:lam2} states that the reconstruction
problem for (\ref{eq:assymetric}) is solvable when $b \lam_2(M) > 1$ and
$\delta_1$ is sufficiently small; Similarly, the reconstruction problem for
(\ref{eq:potts}) is solvable when $b \lam_2(M) > 1$ and $q$ is sufficiently
large.
In this section we improve the existing bounds for non-reconstruction
by showing that:

\begin{Proposition} \label{prop:assymetric_no_reconstruct}
Let $M$ be defined by the transition matrix (\ref{eq:assymetric}).
Then the reconstruction problem for $T_b$ is unsolvable when
\begin{equation} \label{eq:assymetric_no_reconstruct}
b \frac{(\delta_2 - \delta_1)^2}{\delta_2 + \delta_1} \leq 1.
\end{equation}
Similarly, the reconstruction problem for a general tree $T$ is unsolvable if
\begin{equation} \label{eq:assymetric_no_reconstruct_general}
\br(T) \frac{(\delta_2 - \delta_1)^2}{\delta_2 + \delta_1} < 1.
\end{equation}

\end{Proposition}

\medskip\noindent
{\bf Proposition \ref{prop:potts_no_reconstruct}.}
{\it
Let $M$ be the channel (\ref{eq:potts}).
Then the reconstruction problem for $T_b$ is unsolvable when
\begin{equation} \label{eq:potts_no_reconstruct}
b \frac{(1 - q \delta)^2}{1 - (q-2) \delta} \leq 1.
\end{equation}
Similarly, the reconstruction problem for a general tree $T$ is unsolvable if
\begin{equation} \label{eq:potts_no_reconstruct_general}
\br(T) \frac{(1 - q \delta)^2}{1 - (q-2) \delta} < 1.
\end{equation}
}

Propositions \ref{prop:potts_no_reconstruct} and
\ref{prop:assymetric_no_reconstruct} follow from the following theorem.

\begin{Theorem} \label{thm:no_reconstruct}
Let $M$ be a channel and let $i,j$ be two states such that for all
$k \notin \{i,j\}$ it holds that $\M_{i,k} = \M_{j,k}$,
and there exist $0 \leq \eps \leq 1$ and
$\beta \geq 0,\gamma \geq 0$ such that
\begin{equation} \label{eq:two_state_branch}
\left(
\begin{array}{ll}
\M_{i,i} & \M_{i,j} \\
\M_{j,i} & \M_{j,j} \end{array}
\right)
= \alpha \left( \begin{array}{ll} 1 - \eps & \eps \\ \eps & 1 - \eps
            \end{array} \right)
  + \left( \begin{array}{ll} \beta & \gamma \\ \beta & \gamma \end{array}
        \right).
\end{equation}
If $T_b$ is the $b$-ary tree and $b \alpha (1 - 2 \eps)^2 \leq 1$,
then $\lim_{n \to \infty} D_V(\P_n^i,\P_n^j) = 0$.
Similarly, if $T$ is a general tree with
$\br(T) \alpha (1 - 2 \eps)^2 < 1$, then
$\lim_{n \to \infty} D_V(\P_n^i,\P_n^j) = 0$.
\end{Theorem}

\prooft{ of Proposition \ref{prop:assymetric_no_reconstruct}}
If $\delta_1 + \delta_2 \leq 1$ write:
\[
\left( \begin{array}{ll} 1 - \delta_1 & \delta_1
                           \\ 1 - \delta_2 & \delta_2
           \end{array} \right)
=
(\delta_1 + \delta_2)
\left( \begin{array}{ll}
\frac{\delta_2}{\delta_1 + \delta_2} & \frac{\delta_1}{\delta_1 + \delta_2}\\
\frac{\delta_1}{\delta_1 + \delta_2} & \frac{\delta_2}{\delta_1 + \delta_2}
\end{array} \right)
+
\left( \begin{array}{ll} 1 - \delta_1  - \delta_2 & 0
                           \\ 1- \delta_1 - \delta_2 & 0
           \end{array} \right)
\]
and the proposition follows from Theorem \ref{thm:no_reconstruct}.
Otherwise, write
\[
\left( \begin{array}{ll} 1 - \delta_1 & \delta_1
                           \\ 1 - \delta_2 & \delta_2
           \end{array} \right)
=
\left( \begin{array}{ll} \delta_1' & 1 - \delta_1'
                           \\ \delta_2' & 1 - \delta_2'
           \end{array} \right)
\]
where $\delta_1' + \delta_2' < 1$, and use a similar decomposition.
$\QED$

\prooft{ of Proposition \ref{prop:potts_no_reconstruct}}
For all $i \neq j$ we have
\[
\left( \begin{array}{ll} \M_{i,i} & \M_{i,j} \\ \M_{j,i} & \M_{j,j}
\end{array} \right)
=
\left( 1 - (q-2) \delta \right)
\left( \begin{array}{ll}
\frac{1 - (q-1) \delta}{1 - (q-2) \delta} & \frac{\delta}{1 - (q-2) \delta}\\
\frac{\delta}{1 - (q-2) \delta} & \frac{1 - (q-1) \delta}{1 - (q-2) \delta}
\end{array} \right).
\]
and the proposition follows from Theorem \ref{thm:no_reconstruct}.
$\QED$

\prooft{of Theorem \ref{thm:no_reconstruct}}
Assume that the label of the root is chosen to be $i$ or $j$ with
probability $1/2$ each and let $X$ denote this random label.
For a set $W$, we denote by $\widetilde{Y}_W$ the random labeling of
the vertices of $W$,
and $\widetilde{Y}_n = \widetilde{Y}_{L_n}$.
We will show that
\begin{equation} \label{eq:mut_inf_zero}
\lim_{n \to \infty} I(X,\widetilde{Y}_n) = 0,
\end{equation}
where $I$ is the mutual information operator, i.e.,
$I(X,Y) = H(X) + H(Y) - H(X,Y)$ where $H$ is the entropy operator
(see \cite{CT} for more background).
Equation (\ref{eq:mut_inf_zero}) is equivalent to
$\lim_{n \to \infty} D_V(\P_n^i,\P_n^j) = 0$, see e.g. \cite{M:lam2}.

We will split the proof into three main steps.

\medskip\noindent
{\bf Step 1 \cite{EKPS}}: If $M$ is the binary symmetric channel,
\begin{equation} \label{eq:binsym}
\M = \left( \begin{array}{ll}
1 - \eps & \eps \\
\eps & 1 - \eps \end{array}
\right),
\end{equation}
then for a general tree $T$ and a set $W$,
\begin{equation} \label{eq:sym_inf_bound}
I(X,\widetilde{Y}_W) \leq \sum_{w \in W} (1 - 2 \eps)^{2 |w|}.
\end{equation}
This is Theorem 1.3 in \cite{EKPS}.

\medskip\noindent
{\bf Step 2:} We will show that it suffices to prove the theorem assuming
that
\begin{equation} \label{eq:two_state_branch_reduction}
\left(
\begin{array}{ll}
\M_{i,i} & \M_{i,j} \\
\M_{j,i} & \M_{j,j} \end{array}
\right)
= \alpha \left( \begin{array}{ll} 1 - \eps & \eps \\ \eps & 1 - \eps
            \end{array} \right)
\end{equation}
instead of (\ref{eq:two_state_branch}).
Indeed, assume that the theorem is true under the condition
(\ref{eq:two_state_branch_reduction}). Now let $M$ be a channel which satisfies
(\ref{eq:two_state_branch}). Consider the following auxiliary channel $N$.
The channel has the state space $\A' = \A \cup \{i^{\ast},j^{\ast}\}$ and the
following transition matrix $\N$:
\begin{enumerate}
\item
For all $\ell \not \in \{i,j\}$ and all $\ell' \in \A$, set
$\N_{\ell,\ell'} = \M_{\ell,\ell'}$.
\item
\[
\left(
\begin{array}{ll}
\N_{i,i} & \N_{i,j} \\
\N_{j,i} & \N_{j,j} \end{array}
\right)
= \alpha \left( \begin{array}{ll} 1 - \eps & \eps \\ \eps & 1 - \eps
            \end{array} \right).
\]
\item
\[
\left(
\begin{array}{ll}
\N_{i,i^{\ast}} & \N_{i,j^{\ast}} \\
\N_{j,i^{\ast}} & \N_{j,j^{\ast}} \end{array}
\right)
= \left( \begin{array}{ll} \beta & \gamma \\ \beta & \gamma
            \end{array} \right).
\]
\item
For all $\ell$, set $\N_{i^{\ast},\ell} = \N_{i,\ell}$ and
$\N_{j^{\ast},\ell} = \N_{j,\ell}$.
\end{enumerate}
It is clear that if the original channel $M$ satisfied the conditions of
Theorem \ref{thm:no_reconstruct}, then the channel $N$ satisfies these
conditions with (\ref{eq:two_state_branch_reduction})
replacing (\ref{eq:two_state_branch}). By our assumption
this implies that $\lim_{n \to \infty} I(X,\widehat{Y}_n) = 0$, where
$\widehat{Y}_n$ is the labeling of level $n$ for $N$.

Writing $\widetilde{Y}_n$ for the labeling of level $n$ for $M$,
we note that $\widetilde{Y}_n$ may be obtained from $\widehat{Y}_n$ by replacing
each occurrence of $i^{\ast}$ by $i$ and each occurrence of $j^{\ast}$ by
$j$. By the Data Processing Lemma (see e.g. \cite{CT}) it now follows
that
\[
I(X,\widetilde{Y}_n) \leq I(X,\widehat{Y}_n) \to 0,
\]
as needed.

\medskip\noindent
{\bf Step 3:} We the theorem assuming that
(\ref{eq:two_state_branch_reduction}) holds.
We begin by introducing two random variables $Z$ and $Y$.
Recall that $\pa(v)$ is the path from $\rho$ to $v$.
We let:
\begin{equation} \label{eq:defZ}
Z = \left\{ (v,\sigma_v) :
    \exists w \in \pa(v), \sigma_w \notin \{i,j\} \right\},
\end{equation}
and given a set $W$, let
\begin{equation} \label{eq:defY}
Y_{W} = \left\{ (v,\sigma_v) : v \in W,
    \forall w \in \pa(v), \sigma_w \in \{i,j\} \right\}.
\end{equation}
We denote $Y_{L_n}$ by $Y_n$.
Roughly speaking, $Y_n$ contains all information on
the $i \to j$ and $j \to i$ process from the root to level $n$;
$Z$ contains all information which is independent of this process.

By the Data Processing Lemma
\begin{equation} \label{eq:YYZ}
I(X,\widetilde{Y}_n) \leq I(X,(Y_n,Z)).
\end{equation}

Since for all $\ell \notin \{i,j\}$ it holds that
$\M_{i,\ell} = \M_{j,\ell}$, it follows that
$Z$ is independent of $X$. We therefore obtain
\begin{eqnarray} \nonumber
I(X,(Y_n,Z))    &=& H(X) + H(Y_n,Z) - H(X,Y_n,Z)
            = H(X | Z) + H(Y_n | Z) - H(X,Y_n | Z) \\ \label{eq:condZ}
        &=& \E_z I(X,Y_n | Z = z).
\end{eqnarray}

We will show that for almost all $z$,
\begin{equation} \label{eq:lim_information}
\lim_{n \to \infty} I(X,Y_n | Z = z) = 0.
\end{equation}
Since $0 \leq I(X,Y_n | Z = z) \leq 1$ for all $z$,
this would imply the theorem by (\ref{eq:condZ}) and (\ref{eq:YYZ}).
By the Data Processing Lemma, in order to prove (\ref{eq:lim_information}),
it suffices to prove that for a.e. $z$ there exist cutsets $W_n$ for which
\begin{equation} \label{eq:lim_information_cutset}
\lim_{n \to \infty} I(X,Y_{W_n} | Z = z) = 0.
\end{equation}
(Recall that $W$ is a cutset if $W$ intersects every infinite path emanating
from $\rho$;
If $W_n$ is a cutest, then
and if for all $v \in W_n$ we have $\ell \leq |v| \leq \ell'$,
then by the Data Processing Lemma it follows
\[
I(X,Y_{\ell'}) \leq I(X,Y_{W_n}) \leq I(X,Y_{\ell}),
\]
so (\ref{eq:lim_information})
is equivalent to (\ref{eq:lim_information_cutset})).

The key observation is noting that given $Z=z$ we have a broadcast process
with the binary symmetric channel (\ref{eq:binsym}) on the tree
\[
T_z = \left\{v : \forall w \in \pa(v), \sigma_w \in \{i,j\} \right\}.
\]
Therefore by Step 1, in order to prove (\ref{eq:lim_information_cutset})
 it suffices to show that for a.e. $z$, there exist cutsets
$W_n$ such that
\begin{equation} \label{eq:expZ}
\lim_{n \to \infty} \sum_{w \in W_n} (1 - 2 \eps)^{2 |w|} = 0.
\end{equation}
However, for a general tree $T$, we have $\br(T_z) = \alpha \br(T)$
for a.e. $z$. So if $\alpha \br(T) (1 - 2 \eps)^2 < 1$, then (\ref{eq:expZ})
hold for a.e. $z$ for appropriate $W_n$'s.
Equation (\ref{eq:expZ}) holds for $T_d$ when
$b \alpha (1 - 2 \eps)^2 = 1$ by \cite[Theorem 3]{LP}.
$\QED$

\section{$\br(T) |\lam_2|^2 > 1$ implies reconstruction} \label{sec:general_count}
The following proposition is a consequence of the results of
\cite{KS}.
\begin{Proposition} \label{prop:dlam2}
Suppose that $b |\lam_2(M)|^2 > 1$, then the reconstruction problem is census
solvable.
\end{Proposition}
In this section we prove an extension of the proposition
to general trees.

\medskip\noindent
{\bf Theorem \ref{thm:blam2}.}
{\it
Let $T$ be an infinite tree and $M$ a channel such that
$\br(T) |\lam_2(M)|^2 > 1$, then the reconstruction problem for $M$ on $T$
is solvable.
}

In order to prove this theorem we are going to bound from below
$D_{\chi^2}(\P^{(c),i}_n,\P^{(c),j}_n)$ where
\begin{equation} \label{eq:chidef}
D_{\chi^2}(P,Q) = \frac{1}{2}
\sqrt{ \sum_{\sigma} \frac{2 (P(\sigma) - Q(\sigma))^2}
    {P(\sigma) + Q(\sigma)}} ,
\end{equation}
is the $\chi^2$ distance.
\begin{Lemma} \label{lem:chi_and_V}
\begin{equation} \label{eq:chi_totvar}
D_{\chi^2}^2 \leq D_V \leq D_{\chi^2}.
\end{equation}
\end{Lemma}

\proofs
The first inequality in (\ref{eq:chi_totvar}) follows when we use the
estimate:
\[
\frac{2 |P(\sigma) - Q(\sigma)|}{P(\sigma) + Q(\sigma)} \leq 2,
\]
while the second follows from Cauchy-Schwartz when we write:
\[
\frac{1}{2} |P(\sigma) - Q(\sigma)| =
\frac{|P(\sigma) - Q(\sigma)|}{\sqrt{2(P(\sigma) + Q(\sigma))}}
\sqrt{\frac{P(\sigma) + Q(\sigma)}{2}}.
\]
$\QED$

We prove Theorem \ref{thm:blam2} by constructing
linear estimators of the root variable
for finite trees and then evaluating the first and second moments
of these estimators.

Abbreviate $\lam = \lam_2(M)$, and take $v$ to be a right eigenvector
(that is, $\M v = \lam v$) with $|v|_2 = 1$.
We use the notion of {\sl flows} and view the tree as an
{\sl electrical network}. We refer the reader to \cite{DS} and \cite{EKPS} for
definitions and more background.
For a finite tree $T$, we consider $T$ as an electrical network,
where we assign the edge $e$ the resistance
\[
R(e) = (1 - |\lam|^2) |\lam|^{-2 |e|}.
\]
For a vertex $x$, we let $c_x = e_i$
if the label of $x$ is $i$, i.e. if $\sigma_x = i$.

We say that a set of vertices $W$ is an {\em anti-chain} if no vertex in $W$ is
a descendant of another.

\begin{Lemma} \label{lem:resist_finite}
Let $T$ be a finite tree, $W$ an anti-chain in $T$, and $\mu$ a unit flow
from $\rho$ to $W$. Consider the estimator
\begin{equation} \label{eq:defS_mu}
S_{\mu} = \sum_{x \in W} \frac{\mu(x) c_x \cdot v}{\lam^{|x|}}.
\end{equation}
Then for all $\ell \in \A$,
\begin{equation} \label{eq:exp_S_mu}
\E^{\ell}[S_{\mu}] = e_{\ell} \cdot v,
\end{equation}
and there exist a constant $0 < c(M) < 1$ which depends on $M$ only such that
\begin{equation} \label{eq:resistence_bound}
c(M) (1 + \Reff(\rho \lra W))
\leq  \min_{\mu} \min_l \E^{\ell}[|S_{\mu}|^2] \leq
      \min_{\mu} \max_l \E^{\ell}[|S_{\mu}|^2] \leq
      (1 + \Reff(\rho \lra W)).
\end{equation}
\end{Lemma}

\proofs
If $x$ is at level $n$, then
\begin{equation} \label{eq:expected_node_count}
\E_n^i [c_x] = e_i \M^n.
\end{equation}
It follows from (\ref{eq:expected_node_count})
that for every column vector $v$ we have
\begin{equation}
\E_n^i[c_x \cdot v] = e_i \M^n v.
\end{equation}
If we take $v$ to be a vector such that $v$ is an eigenvalue of $\M$
which corresponds to $\lam = \lam_2(M)$, we obtain,
\begin{equation}
\E_n^{\ell}[c_x \cdot v] = \lam^n e_{\ell} \cdot v.
\end{equation}
Now (\ref{eq:exp_S_mu}) follows by linearity.

We let $\pa(x)$ be the path from $\rho$ to $x$ (more generally let
$\pa(x,y)$ be the path from $x$ to $y$).
We let $x \wedge y$, be the {\sl meeting point} of $x$ and $y$, that is, the
vertex farthest from the root $\rho$ on $\pa(x) \cap \pa(y)$.
We have:
\begin{eqnarray} \label{eq:S_corr}
\E^{\ell}[|(c_x \cdot v)  (c_y \cdot v)|] &=& \nonumber
\sum_i
\P^{\ell} [c_{x \wedge y} = e_i] \E[|(c_x \cdot v) (c_y \cdot v)| \, | c_{x \wedge y} = e_i] \\ &=&
\sum_i \P^{\ell}[c_{x \wedge y} = e_i] |e_i \cdot v|^2 |\lam|^{d(x,y)} \leq |\lam|^{d(x,y)}.
\end{eqnarray}
Therefore,
\begin{equation} \label{eq:S_2nd_moment}
\E^{\ell}[|S_{\mu}|^2] \leq
\sum_{x,y} \frac{\mu(x) \mu(y)}{|\lam|^{|x|} |\lam|^{|y|}}
\E^{\ell}[|(c_x \cdot v)(c_y \cdot v)|]
\leq \sum_{x,y} \frac{\mu(x) \mu(y)}{|\lam|^{2 |x \wedge y|}}.
\end{equation}
Since
\[
\frac{1}{|\lam|^{2 |u|}} = 1 + \sum_{e \in \path(u)} R(e),
\]
if follows that for all $\ell$,
\begin{equation} \label{eq:S_mu_2nd}
\E^{\ell}[|S_{\mu}|^2]
\leq (1 + \sum_e R(e) \sum_{x,y \in W}
\one_{e \in \path(x \wedge y)} \mu(x) \mu(y))
= (1 + \sum_e R(e) \mu^2(e)).
\end{equation}
When we take the minimum in (\ref{eq:S_mu_2nd}) we obtain,
\begin{equation} \label{eq:resistence_bound1}
\min_{\mu} \max_{\ell}
\E^{\ell}[|S_{\mu}|^2] \leq \min_{\mu} (1 + \sum_e R(e) \mu^2(e))
= \left(1 + \Reff(\rho \lra W) \right),
\end{equation}
where the equality in (\ref{eq:resistence_bound1}) is Thompson's principle
(see Doyle and Snell \cite{DS}).
Equation (\ref{eq:resistence_bound1}) is the left hand side inequality in
(\ref{eq:resistence_bound}). The proof of the other inequality in
(\ref{eq:resistence_bound}) follows in a similar way.
$\QED$

\prooft{of Theorem \ref{thm:blam2}}
Since $|\lam|^2 \br(T) > 1$ it follows that
\begin{equation} \label{eq:ask_yuval}
\Reff(\rho \lra \infty) = \sup_n \Reff(\rho \lra L_n) < \infty.
\end{equation}
We consider linear estimators as in Lemma \ref{lem:resist_finite} and note
that if $v$ is an eigenvector as in Lemma \ref{lem:resist_finite}, then
\begin{equation} \label{eq:def_tilde_C}
\widetilde{C}(M) = \max_{i,j} | e_i \cdot v - e_j \cdot v | > 0.
\end{equation}
We fix a level $n$, take $i$ and $j$ such that we obtain the maximum in
(\ref{eq:def_tilde_C}),
and consider a linear estimator as in Lemma \ref{lem:resist_finite} for
$W = L_n$, such that
the upper bound in (\ref{eq:resistence_bound}) holds.
We obtain that
\begin{equation} \label{eq:lower_chi}
\left| \sum_{\sigma} S_{\mu}(\sigma)
(\P_n^i[\sigma] - \P_n^j[\sigma]) \right| \geq
\widetilde{C}(M).
\end{equation}
On the other hand, by Cauchy-Schwartz we obtain that
\begin{eqnarray} \nonumber
\left(  \sum_{\sigma} S_{\mu}(\sigma)
(\P_n^i[\sigma] - \P_n^j[\sigma]) \right)^2 &\leq&
\sum_{\sigma} |S_{\mu}(\sigma)|^2 (\P_n^i[\sigma] + \P_n^j[\sigma])
\sum_{\sigma} \frac{(\P_n^i[\sigma] - \P_n^j[\sigma])^2}
              {\P_n^i[\sigma] + \P_n^j[\sigma]} \\ \label{eq:upper_chi}
&\leq& 4 (1 + \Reff(\rho \lra L_n)) D^2_{\chi^2}(\P_n^i,\P_n^j).
\end{eqnarray}
Combining (\ref{eq:lower_chi}) and (\ref{eq:upper_chi}) we obtain that
for all $n$
\begin{equation} \label{chi_estimate}
D^2_{\chi^2}(\P_n^i,\P_n^j) \geq
\frac{\widetilde{C}(M)^2}{4 (1 + \Reff(\rho \lra L_n))}.
\end{equation}
So the theorem follows by Lemma \ref{lem:chi_and_V}.
\QED

\prooft{of Proposition \ref{prop:dlam2}}
For the $b$-ary tree $T_b$ we have $\br(T_b) = b$.
Moreover, the flow $\mu$ which minimizes $\sum_e R(e) \mu^2(e)$ in
(\ref{eq:resistence_bound1}) is the flow which satisfies
$\mu(e) = b^{-|e|}$. This implies that we may take $S_n$ to be
\[
S_n = \frac{c_n \cdot v}{d^n \lam^n}.
\]
In particular $S_n$ is a function of $c_n$.
Now, arguing as in the proof of Theorem \ref{thm:blam2} we see that there
exist $i,j$ for which $D_V(\P_n^{(c),i},\P_n^{(c),j})$ is bounded away from
$0$ for all $n$.
$\QED$

\section{ Census reconstruction fails when $b |\lam_2|^2 < 1$}
\label{sec:easy_count}

In this section we prove non-reconstruction when $b |\lam_2|^2 < 1$.
The proof relies on the Kesten-Stigum theorem \cite{KS}.
The following is an immediate consequence of Theorem 2.3 of \cite{KS}.

\medskip\noindent
{\bf Kesten-Stigum CLT \cite{KS}}:
{\em
Let $M$ be a transition matrix such that $b |\lam_2(M)|^2 < 1$ and
let $\pi$ be the stationary distribution for $M$, i.e., the normalized
left eigenvector for the eigenvalue 1. Then for any
vector $v$ of $M$ which is orthogonal to $\pi$, and for all $j \in \A$,
\begin{equation} \label{eq:normal_law}
\frac{(c^j_n,v)}{b^{n/2}} \to \NC(0,\sigma),
\end{equation}
where $\NC$ denotes a normal random variable and $\sigma$ does not depend on $j$.
Moreover, the convergence rate may be bounded in term of $|v|_2$ only.
}

An immediate consequence is:
\begin{Proposition} \label{prop:normal_law}
For all $i,j$ and $\eps > 0$, one may couple $c^j_n$ and $c^i_n$
in such a way that,
\begin{equation} \label{eq:normal_law2}
\lim_{n \to \infty} \P[\|c^j_n - c^i_n \| > \eps b^{n/2}] = 0.
\end{equation}
\end{Proposition}

\begin{proof}
If follows from the Kesten Stigum theorem that
\[
\frac{c^{\ell}_n- \pi b^n}{b^{n/2}} \to \NC
\]
where $\NC$ is a normal variable which does not depend on
$\ell$.
The claim follows, as for all $\delta > 0$ we can couple both $c^i_n$
and $c^j_n$ with the same normal variable $b^{n/2} \NC$ in such a way that
\[
\P[\| c^{\ell}_n - b^{n/2} \NC \| > \eps b^{n/2}] < \delta,
\]
for both $\ell=i$ and $\ell=j$.
\end{proof}

\medskip\noindent
{\bf Theorem \ref{thm:non_count_rec}.}
{\it
If $b |\lam_2(M)|^2 < 1$, then the reconstruction problem is not
census solvable for the $b$-ary tree and the channel $M$.
}

The proof is easier when all the entries of $\M$ are strictly positive,
in which case the following lemma is trivial.
\begin{Lemma}
Let $M$ be irreducible and aperiodic, then there exists an $h$ such that
for all $i$ and $j$ all $n \geq h$ and all $v$,
\begin{equation} \label{eq:rw_steps}
\P[c^i_n = v] > 0 \mbox{ iff }\P[c^j_n = v] > 0.
\end{equation}
\end{Lemma}

\begin{proof}
Since $(\M_{i,j})_{i,j=1}^k$ is irreducible and aperiodic it follows that there exists an $\ell$
such that the matrix $\M^\ell$ is strictly positive. Thus, for every two states $i$ and $j$ we may construct
a $\ell$ level tree such that the root of the tree is labeled by $i$ and all of the leaves are labeled by $j$.
Let $h = \ell(k + 4)$ and $n \geq h$.
Suppose that $\P[c^i_n = v] > 0$. Then there exists a labeling of the tree of $n$ levels
which has $i$ at the root and $\sigma$ at level $n$, and such the census of $\sigma$
is $v$. We will prove the lemma by constructing a labeling where the root is labeled by $j$,
level $n$ is labeled by $\tau$, and the census of $\tau$ is $v$.

We denote by $x_1, \ldots, x_{b^{n - \ell}}$ the vertices at level $n - \ell$.
For $x_t$ we define $c(x_t)$ to be the census of the subtree rooted at $x_t$ for the labeling $\sigma$.
Let $u=(u_1,\ldots,u_k)$ be a non-negative vector which satisfies $\sum u_i = b^{\ell}$.
Define $c(\sigma,u) = |\{t : c(x_t) = u\}|$. Note that if $c(\sigma,u) > 0$ then there exists a labeling
of the $\ell$ levels tree, denoted $\tau(u)$, which has $u$ as its census and such that the label of the
root for this labeling is $\sigma_{\rho}(u)$.

Since
\[
v - \sum_u \left( c(\sigma,u) \mod b^{\ell} \right) u = 0 \mod b^{\ell},
\]
it follows that
\[
v = b^{\ell} w + \sum_u \left( c(\sigma,u) \mod b^{\ell} \right) u,
\]
for some integer valued vector $w$.

Note that if for root value $j$, we could label $\sum_u \left( c(\sigma,u) \mod b^{\ell} \right)$ of the vertices of
level $n - \ell$ by labels $\sigma_{\rho}(u)$ with multiplicity $c(\sigma,u) \mod b^{\ell}$, then we are done.
Since using these vertices, it is possible to build the $\sum_u \left( c(\sigma,u) \mod b^{\ell} \right) u$ part of the
census. Now using the other vertices at level $n - \ell$ it is possible to build the $b^{\ell} w$ part of the census
(whatever labels these vertices have).
Note that
\[
\sum_u \left( c(\sigma,u) \mod b^{\ell} \right) \leq b^{\ell} |\{ u : \forall i, u_i \geq 0, \sum_{i=1}^k u_i = b^{\ell}\}|
\leq b^{\ell} b^{k \ell} = b^{(k+1) \ell}.
\]
Therefore by assigning $b^{(k+1) \ell}$ of the vertices at level $n - 2 \ell$ the task of producing the
prescribed labels at level $n - \ell$ we obtain the required result.
\end{proof}

\prooft{of Theorem \ref{thm:non_count_rec}}
We take $h$ such that (\ref{eq:rw_steps}) holds for $n \geq h$.
We may write:
\begin{equation} \label{eq:dec_CLT_LCLT}
c^j_{n+h} = \sum_{i=1}^k S^i(c^j_n(i)),
\end{equation}
where $S^i$ is a random walk on ${\Z^k}$ which satisfy:
\[
\P[S^i(t+1) = S^i(t) + v] = \P[c^i_h = v].
\]
By (\ref{eq:rw_steps}), all the $S^i$ are random walks with the same support.

By the local CLT if follows that if $|c^j_n - c^i_n| < \eps \sqrt{b^n}$,
it is possible to couple $c^j_{n+h}$ and $c^i_{n+h}$ with probability $f(\eps)$
where $f(\eps) \to 0$ as $\eps \to 0$. Now the result follows from Proposition \ref{prop:normal_law}.

$\QED$

\section{Census reconstruction for Potts models fails at criticality}
 \label{sec:potts_count}

In section \ref{sec:easy_count} we presented a proof that $b |\lam_2|^2 < 1$
implies that the reconstruction problem is not census-solvable.
We believe that this result is also valid when $b |\lam_2|^2 = 1$.
In this section we prove that this is the case
when $M$ is the $q$ state Potts model or the Ising model
with external field.


\medskip\noindent
{\bf Theorem \ref{thm:no_count_reconstruct}.}
{\it
Let $M$ be the $q$ state Potts model, or the asymmetric Ising model,
and suppose that $b \lam_2^2(M) \leq 1$, then the reconstruction problem
is unsolvable.
}

\proofs
We prove the theorem for the $3$-state Potts model, the general proof being similar.
It is convenient to denote the states of the channel
by $-1,0$ and $1$.
We show that for all $i,j \in \{-1,0,1\}$ there exist a
coupling such that $\P[c^i_n = c^j_n] \to 1$ as
$n \to \infty$.
By symmetry, it suffices to find such coupling for $c^1_n$ and
$c^{-1}_n$. We denote by $\sigma^{+}$ the coloring with $1$ at the root,
and $\sigma^{-}$ for the coloring with $-1$ at the root.

During all the steps of the coupling we require that for all $v$ we have $\sigma^{-}_v = 0$ iff
$\sigma^{+}_v = 0$.
The main step in proving the existence of the required coupling is given in
the following lemma.
\begin{Lemma} \label{lem:iter_coupling}
There exits $p^{\ast} > 0$ such that:
Given $(\sigma_v)_{|v| \leq n}$ and $(\tau_v)_{|v| \leq n}$
which are coupled in such a way that $\sigma_v = 0$ iff $\tau_v = 0$,
there exits $N \geq n$ and a coupling procedure for $(\sigma_v)_{|v| \leq N}$ and
$(\tau_v)_{|v| \leq N}$ such that $\sigma_v = 0$ iff $\tau_v = 0$ and
such that $\P[c^{\sigma}_N = c^{\tau}_N] \geq p^{\ast}$ (where $c^{\sigma}_N$ and $c^{\tau}_N$ are the
$\sigma$ and $\tau$ level $N$ census vectors respectively).
\end{Lemma}
By iterating Lemma \ref{lem:iter_coupling}, the coupling probability $\P[c^1_n = c^{-1}_n]$
may be made arbitrarily close to $1$ for large $n$.
$\QED$

\prooft{of Lemma \ref{lem:iter_coupling}}
At level $n$ there are two types of vertices: those for which $\sigma_v = \tau_v$, and
those for which $\sigma_v = -\tau_v \in \{-1,1\}$.
For all the vertices $w$ whose ancestors $v$ at level $n$ satisfy $\sigma_v = \tau_v$,
we let $\sigma_w = \tau_w$.

Denote the vertices $v$ at level $n$ which satisfy
$\sigma_v = -\tau_v \in \{-1,1\}$ by $v_1,\ldots,v_{\ell}$.
Let $T$ be the graph which consists of the subtrees rooted at $v_1,\ldots,v_{\ell}$.

We think of $T$ as drawn in the plan,
and denote the vertices of level $m$ by $v^m_1,\ldots,v^m_{\ell b^{n-m}}$.
We will slightly abuse the notation by redefining
\begin{equation} \label{eq:redef}
\begin{array}{ll}
c^{\sigma}_m(i) = |\{j : 1 \leq j \leq \ell b^{n-m}, \sigma(v^m_j) = i\}|, &
c^{\tau}_m(i) = |\{j : 1 \leq j \leq \ell b^{n-m}, \tau(v^m_j) = i\}|.
\end{array}
\end{equation}
Note that $c^{\sigma}_m = c^{\tau}_m$ in the old definition if and only if it does in the new one.

We define $X^{\sigma}_m = c^{\sigma}_m(1) - c^{\sigma}_m(-1)$ and
$X^{\tau}_m = c^{\tau}_m(1) - c^{\tau}_m(-1)$. Note that if a coupling of $\sigma$
and $\tau$ satisfies $\sigma_v = 0$ iff $\tau_v = 0$, then we have
$c^{\sigma}_m = c^{\tau}_m$ if and only if $X^{\sigma}_m = X^{\tau}_m$.

We define the coloring of $T$ as a dynamic process starting
at the root, running level by level from left ($v^m_1$) to right ($v^m_{\ell b^{n-m}}$).
Writing $v'$ for the parent of $v$, we let
\begin{equation} \label{eq:def_mart}
S^{\sigma}(v^m_i) = (b \theta)^{-n} X^{\sigma}_n + \sum_{n \leq k < m}
(b \theta)^{-k} \sum_{j \leq \ell b^{k-n}} (\sigma(v^k_j) - \theta \sigma({v^k_j}')) +
(b \theta)^{-m} \sum_{j \leq i} (\sigma(v^m_j) - \theta \sigma({v^m_j}')),
\end{equation}
where we formally define $\sigma_{\rho'} = 0$.
We define $S^{\tau}$ similarly.
Note that all the terms in (\ref{eq:def_mart}) except the first one, have mean zero,
and therefore both $S^{\sigma}$ and $S^{\tau}$ are martingales.
Also,
\begin{equation} \label{eq:SX}
\begin{array}{ll}
S^{\sigma}(v^m_{\ell b^{m-n}}) = (b \theta)^{-m} X^{\sigma}_m, &
S^{\tau}(v^m_{\ell b^{m-n}}) = (b \theta)^{-m} X^{\tau}_m.
\end{array}
\end{equation}

The coupling consists of a few steps, during all of which we require that for all $v$
we have $\sigma_v = 0$ iff
$\tau_v = 0$.
\begin{enumerate}
\item
{\bf Reflection until we reach $N$ such that:}
\begin{equation} \label{eq:reflect_cond}
-\sqrt{b^N} C_1 \leq X^{\sigma}_N = -X^{\tau}_N \leq C_1 \sqrt{b^N}.
\end{equation}
We label the vertices level by level, left to right using the rule $\sigma_v = - \tau_v$,
until we reach $v^N_i$ with
$|S^{\sigma}(v^N_i)| \leq (1 + \theta) (b \theta)^{-N}$ and $N$ large.
From lemma \ref{lem:mart_inf_sup} below, it follows that such a vertex exits a.s.
Now $X^{\sigma}_N - (b \theta)^N S^{\sigma}(v^N_i)$
is a sum of $\ell b^{N-n} - i$ independent variables with values $-1,0,1$ and expected value $0$.
It therefore follows by the CLT that with probability going to $1$ as $C_1 \to \infty$ we have
$|X^{\sigma}_N| \leq C_1 \sqrt{b^N}$.
Therefore, by continuing the reflection $\sigma(v^N_j) = -\tau(v^N_j)$ for
$i < j \leq \ell b^{N-n}$, we obtain (\ref{eq:reflect_cond}).
We let $Y^{\sigma}_{N+1}(i,j)$ be the number of vertices in the $\sigma$ labeling
which are at level $N+1$, of type $j$ and have a type $i$ parent.
All these variables are binomial and
\begin{equation} \label{eq:expY}
\E[Y^{\sigma , \tau}_{N+1}(i,j) | c^{\sigma , \tau}_N] =
\left( b (1 - \theta)/3 + \delta_{i,j} b \theta \right) c^{\sigma , \tau}_N(i).
\end{equation}
\item
{\bf Labeling the $0$'s of level $N+1$:}
We label all vertices of level $N+1$ with $\sigma(v^{N+1}_j) = \sigma(v^{N+1}_j) = 0$ independently
with probability
\[
\P[\sigma(v^{N+1}_j) = 0 | \sigma((v^{N+1}_j)')] = \P[\tau(v^{N+1}_j) = 0 | \tau((v^{N+1}_j)')].
\]

We thus achieve that $c^{\sigma}_{N+1}(0) = c^{\tau}_{N+1}(0)$.

By the CLT it follows that with
probability going to $1$ as $C_3 \to \infty$,
\[
|Y^{\sigma , \tau}_{N+1}(i,0) - \E[Y^{\sigma , \tau}_{N+1}(i,0) | c^{\sigma,\tau}_N]| \leq C_3 \sqrt{b^{N+1}},
\]
and therefore for $j \neq 0$,
\begin{equation} \label{eq:dev1}
|\E[Y^{\sigma , \tau}_{N+1}(i,j) | c^{\sigma , \tau}_N, Y^{\sigma , \tau}_{N+1}(i,0)] -
 \E[Y^{\sigma , \tau}_{N+1}(i,j) | c^{\sigma , \tau}_N]| \leq C_3 \sqrt{b^{N+1}}.
\end{equation}

\item
{\bf Labeling some of the $\pm$ of level $N+1$.}
Let $i_0$ be the maximizer of $c^{\sigma}_N(i)$.
For $i \neq i_0$ we label all vertices $v$ at level $N+1$ which have a parent of type $i$ and satisfy
$\sigma(v) \neq 0$ by $\pm 1$ using the reflecting coupling $\sigma(v) = -\tau(v)$.
By the CLT it follows that with probability going to $1$ as $C_4 \to \infty$ we have for $i \neq i_0$ and all $j$
\begin{equation} \label{eq:dev2}
|Y^{\sigma , \tau}_{N+1}(i,j) -
\E[Y^{\sigma , \tau}_{N+1}(i,j) | c^{\sigma , \tau}_N, Y^{\sigma , \tau}_{N+1}(i,0)|
\leq C_4 \sqrt{b^{N+1}}.
\end{equation}

\item
{\bf Local CLT for the final coupling.}
We now want to label the remaining vertices with $\pm 1$ to achieve
$X^{\sigma}_{N+1} = X^{\tau}_{N+1}$.
We note that $Z^{\sigma , \tau} = Y^{\sigma , \tau}(i_0,1) - Y^{\sigma , \tau}(i_0,-1)$
are sum of i.i.d. $\pm 1$ random variables and therefore
satisfy the local CLT.
Moreover, by (\ref{eq:dev1}),
it follows that with high probability each of the $Z$'s is a sum of at least
$\ell b^{N+1-n} (1 - \theta)/6$ such variables. In order that $X^{\sigma}_{N+1} = X^{\tau}_{N+1}$ we need that
$\sum_i Y^{\sigma}(i,1) - Y^{\sigma}(i,-1) = \sum_i Y^{\tau}(i,1) - Y^{\tau}(i,-1)$.
By (\ref{eq:expY}),(\ref{eq:dev1}) and
(\ref{eq:dev2}) it follows that for
$i \neq i_0$ and with probability going to $1$ as $C_5 \to \infty$,
\begin{equation} \label{eq:dev3}
| (Y^{\sigma , \tau}(i,1) - Y^{\sigma , \tau}(i,-1)) - b \theta c^{\sigma, \tau}_N(i) |
\leq C_5 \sqrt{b^{N+1}}.
\end{equation}
Similarly, by (\ref{eq:expY}) and (\ref{eq:dev1})
\begin{equation} \label{eq:dev4}
| \E[Y^{\sigma , \tau}(i_0,1) - Y^{\sigma , \tau}(i_0,-1) |
c^{\sigma , \tau}_N, Y^{\sigma , \tau}_{N+1}(i,0)] - b \theta c^{\sigma,\tau}_N(i_0) |
\leq C_5 \sqrt{b^{N+1}}
\end{equation}
Now by (\ref{eq:dev3}),(\ref{eq:dev4}),(\ref{eq:reflect_cond}) and the fact that until now we used reflection,
it follows that
$Z^{\sigma} + W^{\sigma} = Z^{\tau} + W^{\tau}$ where $W^{\sigma , \tau}$
are random variables which satisfies $W^{\sigma} = -W^{\tau}$
and with probability going to $1$ as $C_6 \to \infty$,
$|W^{\sigma} - W^{\tau} - (\E[Z^{\sigma}] - \E[Z^{\tau}])| \leq C_6 \sqrt{b^{N+1}}$.
By the local CLT, it follows the with probability $p^{\ast} > 0$ we may couple $Z^{\sigma}$ and
$Z^{\tau}$ in order to achieve
$X^{\sigma}_N = X^{\tau}_N$.
\end{enumerate}

\begin{Lemma} \label{lem:mart_inf_sup}
A.s. $S^{\sigma}$ changes signs infinitely often.
\end{Lemma}

\proofs
We prove that $X^{\sigma}_n$ changes sign infinitely often a.s.
By the Borel-Cantelli lemma
it suffices to prove that for all $n$ and all
$\sigma(v)_{|v| \leq n}$ there exists $m \geq n$ such that
$ \P[X^{\sigma}_m > 0 | \sigma(v)_{|v| \leq n}] \geq 1/4$.
However, by the Kesten-Stigum theorem \cite{KS}, there exits some $a_m \to \infty$
such that given $\sigma(v)_{|v| \leq n}$,
$X_m/a_m$ converges to a non-degenerate normal random variable with expected value $0$.
So when $m$ is large, $ \P[X^{\sigma}_m > 0 | \sigma(v)_{|v| \leq n}] \geq 1/4$ as needed.
$\QED$

\section{Unsolved problems} \label{sec:open}
\subsection{Critical values}
Except for symmetric binary channels there are no interesting
families of channels for which the critical value for reconstruction
is known.
\begin{Problem}
For the $3$ symbols Potts model (\ref{eq:potts}) find the critical
value for reconstruction.
\end{Problem}
The same question applies to asymmetric binary channels, and to
proper colorings (which correspond to the zero temperature
anti-ferromagnetic Potts model).
Proper colorings of trees were studied in \cite{BW1};
The corresponding transition matrix is
\begin{equation} \label{eq:color}
\M =
 \left(
\begin{array}{lllll} 0 & (q-1)^{-1} & (q-1)^{-1} & \ldots     & (q-1)^{-1}   \\
                    (q-1)^{-1} & 0 &  (q-1)^{-1}       & \ldots \\
                    \vdots & \ldots & \ddots & \vdots            \\
                    (q-1)^{-1} & \ldots & (q-1)^{-1}     & 0
\end{array}
        \right).
\end{equation}
A simple coupling argument (see \cite{BW1}) shows that if $b \leq q-1$,
then the reconstruction problem is unsolvable for colorings of $T_b$.
On the other hand, applying standard coupon-collector
estimates recursively (similarly to Theorem \ref{thm:dist_comb}), it is
easy to see that if $b \geq (1+\epsilon) q \log q$ and $q$ is large, then the
reconstruction problem is solvable. By Proposition \ref{prop:dlam2}
and Theorem \ref{thm:non_count_rec}, the
census reconstruction problem is solvable if $b>(q-1)^2$,
and unsolvable if  $b<(q-1)^2$.

\begin{Problem}
For fixed $q$, for which $b$ is the reconstruction problem solvable
for the channel (\ref{eq:color})?
\end{Problem}

\subsection{Soft inputs - Robust phase transitions}
We like to mention briefly the notion of "robust" phase transition which
first appeared in \cite{PS}.
Consider the usual reconstruction problem, but suppose that the data at the
boundary is given with some additional noise.
The proofs that if $b \lam_2^2(M) > 1$ the reconstruction problem is
(census) solvable are immune to this noise.
However, this may not be the case for the reconstruction problem.
Indeed, we suspect that adding this additional noise
(assuming it is fixed but sufficiently strong) will shift the
phase transition to the point $b \lam_2^2(M) = 1$.
A similar phenomena was proven in \cite{PS} for the phase transition
of uniqueness.
For $n$ and $m$, we denote by $\sigma_{n,m}$ the configuration which
is obtained from $\sigma_n$ by applying the random function
$M^m$ independently
on each of the symbols in $\sigma_n$. We denote by $\P^{\ell}_{n,m}$
the conditional distribution on $\sigma_{n,m}$ given that
$\sigma_{\rho} = \ell$. We then
\begin{Conjecture}
For all $M$ and $b$, such that $b \lam_2^2(M) < 1$, there exists $m$
such that
\[
\sup_{i,j} \lim_{n \to \infty} D_V(\P^i_{n,m},\P^j_{n,m}) = 0.
\]
\end{Conjecture}

\subsection{Monotonicity}
For Potts models (\ref{eq:potts1}), it is easy to see that
if the reconstruction problem is solvable for $q$ and $\eps$
and $q' < q$, then the reconstruction problem is also solvable for $q'$
and $\eps$.
We expect that for fixed $\lam = \lam_2(M)$,
reconstruction is easier when $q$ is larger.
\begin{Conjecture}
Consider two symmetric channels $M_1$ and $M_2$
(as in  (\ref{eq:potts1})) on $q_1$ and $q_2$
symbols respectively, where $q_1 < q_2$.
If $\lam_2(M_1) = \lam_2(M_2)$ and the reconstruction problem
 is solvable for $M_1$, then it is also solvable for $M_2$.
\end{Conjecture}
This is obvious when $q_2$ is a multiple of
$q_1$. Using the reconstruction criteria for the binary symmetric channel on
$2$ symbols, it is easy to prove the conjecture when $q_1 = 2$.




\medskip\noindent
{\bf Acknowledgments:} We thank Olle H\"{a}ggstr\"{o}m, Claire Kenyon,
L\`aszl\`o Lov\`asz, Jeff Steif and Peter Winkler for helpful discussions.


\vspace{.5in}
\noindent {\sc Elchanan Mossel} \newline
Microsoft research \newline
1 Microsoft way \newline
Redmond WA 98584
\newline
 {\tt mossel@microsoft.com}

\vspace{0.5in}
\noindent {\sc Yuval Peres} \newline
Department of Statistics \newline
367 Evans Hall, University of California \newline
Berkeley, CA 94720-3860 \newline
{\tt peres@stat.berkeley.edu}

\begin{thebibliography}{7}

\bibitem{AN} Athreya, K. B. and Ney, P. E. (1972) {\em Branching Processes},
Springer-Verlag.

\bibitem{BW1} Brightwell, G. and Winkler, P. (2001). Random colorings of a Cayley tree, {\sl submitted}.
\bibitem{BW2} Brightwell, G. and Winkler, P. (2000). Gibbs measures and dismantlable graphs,
{\em J. Comb. Theory (Series B)} {\bf 78}, 141--169.
\bibitem{BW3} Brightwell, G. and Winkler, P. (1999). Graph homomorphisms and phase transitions,
{\em J. Comb. Theory (Series B)}, {\bf 77}, 415--435.

\bibitem{BRZ} Bleher, P. M., Ruiz, J. and Zagrebnov V. A. (1995)
On the purity of limiting Gibbs state for the Ising model on the Bethe
lattice, {\sl J. Stat. Phys} {\bf 79}, 473--482.

\bibitem{Ca} Cavender, J. (1978). Taxonomy with confidence.
{\em Math.\  BioSci.\  } {\bf 40},  271--280.

\bibitem{CT} Cover, T. M. and Thomas, J. A. (1991)
{\em Elements of Information Theory}, John Wiley and Sons.

\bibitem{DS}
Doyle, P. G. and Snell, E. J.(1984).
{\em Random walks and Electrical Networks.}
 Carus Math.\ Monographs {\bf 22}, Math.\ Assoc.\ Amer., Washington, D. C.

\bibitem{DZ}
Dembo, A. and Zeitouni O. (1997)
{\it Large Deviations, Techniques and Applications}, Springer.

\bibitem{EKPS} Evans, W., Kenyon, C., Peres, Y. and Schulman L. J. (2000)
Broadcasting on trees and the Ising Model, {\it Ann. Appl. Prob.},
{\bf 10 no. 2}, 410--433.

\bibitem{Fi}  Fitch, W. M. (1971). Toward defining the course of evolution:
 minimum change for a specific tree topology.
{\em Syst.\ Zool.} {\bf 20}, 406--416.

\bibitem{HW} Hajek, B. and Weller, T. (1991). On the maximum tolerable noise
for reliable computation by formulas.
{\em IEEE Trans.\ on Information Theory\/} {\bf 37}(2), 388--391.

\bibitem{Hi} Higuchi, Y. (1977).  Remarks on the limiting Gibbs state
on a (d+1)-tree. {\em Publ.\  RIMS Kyoto Univ.} {\bf 13}, 335--348.

\bibitem{I1} Ioffe, D. (1996a).  A note on the extremality of the disordered
state for the Ising model on the Bethe lattice.
{\em Lett.\  Math.\  Phys.} {\bf 37}, 137--143.

\bibitem{I2}  Ioffe, D. (1996b).  A note on the extremality of the disordered
state for the Ising model on the Bethe lattice. In
{\em Trees}, B. Chauvin, S. Cohen, A. Roualt (Editor).

\bibitem{KMP} Kenyon, C., Mossel, E. and Peres, Y. (2001).
Glauber dynamics on trees and hyperbolic graphs, {\em Preprint}.

\bibitem{KS} Kesten, H. and Stigum, B. P. (1966)
Additional limit theorem for indecomposable multidimensional Galton-Watson
processes, {\it Ann. Math. Statist.} {\bf 37}, 1463--1481.

\bibitem{Lyons} Lyons, R. (1990)
Random walks and percolation on trees. {\em Ann. Probab.} {\bf 18},
931--958.

\bibitem{LP} Lyons, R. and Pemantle R. (1992)
Random walk in a random environment and first-passage percolation on trees.
{\em Ann. Probab.} {\bf 20,1} 125--136.

\bibitem{LP2} Lov\`asz, L. and Winkler P. (1998)
Mixing times, in {\em DIMACS Series in Discrete Mathematics and Theoretical Computer Science}
{\bf 41}, 85--133.

\bibitem{M:recur} Mossel, E. (1998)
Recursive reconstruction on periodic trees,
            {\it Random Structures and algorithms} {\bf 13,1} 81--97

\bibitem{M:lam2} Mossel, E. (2001)
Reconstruction on trees: Beating the second eigenvalue,
           {\it Ann. Appl. Probab.}, to appear.

\bibitem{PS} Pemantle, R. and Steif, J. E. (1999).
Robust phase transitions for Heisenberg and other models on general trees.
{\em Ann. Probab.} {\bf 27 no. 2}, 876--912

\bibitem{Sh} Shamir, A. (1979). How to share a secret?
Communications of the ACM {\bf 22} , 612--613.

\bibitem{Sp} Spitzer, F. (1975). Markov random fields on an infinite tree.
{\em Ann.\ Probab.\  } {\bf 3}, 387--394.

\bibitem{St} Steel, M. (1989).
 Distributions in  bicolored evolutionary trees. {\em Ph.D. Thesis},
 Massey University, Palmerston North, New Zealand.


\end{thebibliography}
\end{document}